\documentclass{scrartcl}

\usepackage{amssymb,amsmath,graphicx,dsfont,xcolor,booktabs,hyperref,units}

\usepackage[bottom]{footmisc}

\usepackage{type1cm}        
\usepackage{newtxtext}       %
\usepackage{newtxmath}       

\usepackage{tikz}
\usepackage{pgfplots}
\newlength \figureheight
\newlength \figurewidth

\pgfplotsset{compat=1.9}
\usetikzlibrary{
  pgfplots.groupplots,
  matrix
}

\newcommand{\vt}{\mathbf{v}}
\newcommand{\xt}{\mathbf{x}}
\newcommand{\Ft}{\mathbf{F}}

\newcommand{\ft}{\mathbf{f}}

\newcommand{\wt}{\mathbf{w}}
\newcommand{\ut}{\mathbf{u}}

\newcommand{\nt}{\vec{n}}

\newcommand{\sigmat}{\boldsymbol{\sigma}}
\newcommand{\Sigmat}{\boldsymbol{\Sigma}}

\renewcommand{\div}{\operatorname{div}}

\newcommand{\FL}{{\cal F}}
\newcommand{\SO}{{\cal S}}

\newtheorem{algo}{Algorithm}

\begin{document}

\title{The candy wrapper problem - a temporal multiscale approach for
  pde/pde systems} 
\author{Thomas Richter~\thanks{~ University of Magdeburg, 
    Universit\"atsplatz 2, 
    39104 Magdeburg, 
    Germany, 
    \texttt{thomas.richter@ovgu.de}
  }
  \and Jeremi Mizerski~\thanks{~ University of Magdeburg, 
    Universit\"atsplatz 2, 
    39104 Magdeburg, 
    Germany, 
    \texttt{jeremi.mizerski@ovgu.de}
  }}
  
\maketitle

\begin{abstract}
  We describe a temporal multiscale approach for the simulation of
  long-term processes with short-term influences involving partial
  differential equations. The specific problem under consideration is
  a growth process in blood vessels. The \emph{Candy Wrapper Process}
  describes a restenosis in a vessel that has previously be widened by
  inserting a stent. The development of a new stenosis takes place on
  a long time horizont (months) while the acting forces are mainly
  given by the pulsating blood flow. We describe a coupled pde model
  and a finite element simulation that is used as basis for our
  multiscale approach, which is based on averaging the long scale
  equation and approximating the fast scale impact by localized
  periodic-in-time problems. Numerical test cases in prototypical 3d
  configurations demonstrate the power of the approach. 
\end{abstract}

\section{Introduction - The candy wrapper problem}

The idea of opening or dilating occluded or narrowed coronary artery
originates in the works of Andreas Gruentzig. First human application
of percutaneous transluminal coronary angioplasty (PTCA) had been
performed on September 16th 1977 at University Hospital in Zurich. The
method was basically just putting the balloon catheter through
narrowing and inflating it \cite{RN36}. The immediate results were
good, only about 1\% of the patients suffered from immediate vessel
closure and myocardial infarct. Later after the interventions 30\% of
the stenosis recurred accompanied by the symptoms of angina of the
intensity close to those from before the intervention. That happened
usually from 30 days to 6 months form the intervention~\cite{RN39}. At
that time the cardiologists were convinced that only about 10\% of all
the patients will be suitable for the method and the rest of coronary
artery disease cases had to be referred to cardiac surgery for by-pass
grafting. 
The remedy for the situation were to be stents intended as an internal
scaffold for the artery to maintain it’s patency. The method was
introduced in 1986 with some success~\cite{RN40}. Soon after that a new set of
complications came into the attention. The early and late onset of
thrombosis started to haunt the patients undergoing procedures of bare
metal stent (BMS) implantation. The BMS coped also with the problem of
intimal hypertrophy which resulted in in-stent stenosis. From that
moment on the era of drug eluting stents (DES) begins. Throughout the
90s different companies try different chemical compounds. The first
successful application was reported by Serruys in 1998~\cite{RN41}. That
however did not solve the problem entirely and resulted in even more
complex set of complications~\cite{RN42,RN45}. The platelet dependent thrombosis
resulted in explosion of anti-platelet drug development in following
years. The problem defined as a ``restenosis of treatment margins'' or
``candy wrapper'' phenomenon was described by radiologists trying to
apply the oncological brachytherapy principles to the neointimal
overgrowth inside BMS~\cite{RN43}. Soon after that the molecular bases of the
process started to be extensively studied~\cite{RN44}. The issue of stent edge
stenosis had not been resolved by introduction of new materials and
coatings~\cite{RN49,RN55}. The biological effects of flow properties have been
studied extensively since the introduction of extracorporeal
circulatory system in early 50s. The body of evidence built on that
experience showed large interdependencies between the local flow
properties and the tissue response. The research areas branched
towards optimization of stent struts geometry~\cite{RN46} and usage of
different cytostatic drugs as a stent coating material~\cite{RN47}. The key
elements of the milieu created by stents are usually considered
separately. Some computational models allow to recreate and integrate
more elements into the system~\cite{RN53,RN57}. By means of computer
simulations the researchers were able to simulate not only fluid
dynamics around the stented area but also the effects of drug
diffusion into the arterial walls~\cite{RN54,RN56}. The edge restenosis
phenomenon however did not find its’ conclusive description. To fully
understand that complex phenomenon we need to take the arterial wall
mechanics and fluid-structure interactions into consideration.
The specific challenge that is tackled in this work is the temporal
multiscale character of this problem: While restenosis occurs after
months, the driving mechanical forces come from the pulsating blood
flow that requires a resolution in the order of centiseconds. Direct
simulations of this long-term process are not feasible and we present
temporal multiscale methods aiming efficient predictions.

\section{Model configuration}

In this section we will briefly describe the mathematical model used
to decribe the stenosis growth effects. Medical, biological and
chemical processes are strongly simplified. They do however still
contain the specific couplings and scales that are characteristic for
the underlying problem. We choose problem parameters as close to the
medical configuration as possible and as known, which is an issue
since good data is difficult to measure and only sparsely available. 

The most important simplification in our present computational model
is the assumption of a rigid vessel wall. Although deformation by
dynamical fluid-structure interactions are small it is well known that
the effects of elasticity should be taken into effect for an
appropriate depiction of wall stresses, which are an essential
ingredient in triggering stenosis growth. However, we give an outlook
on techniques that are suitable to substantially increase the
efficiency in medical fluid-structure interaction simulation that
suffer from special instabilities by the added-mass effect due to
similar masses of fluid and solid~\cite{CausinGerbeauNobile2005}. 

\subsection{Governing equations}

We consider a system of partial differential equations that is
inspired
by~\cite{YangJaegerNeussRaduRichter2015,YangRichterJaegerNeussRadu2017},
where a complex system of partial differential equations describing
the interaction of mechanical fluid-structure interactions with
bio/chemical reactions and active growth and material deformation is
introduced. The 
mechanical system is described by a nonlinear fluid-structure
interaction model, where the blood is modeled as incompressible
Newtonian fluid, which is an adequate choice for the vessel sizes under
consideration
\begin{equation}\label{fluidproblem}
  \rho_f \big(\partial_t\vt + (\vt\cdot\nabla)\vt\big) -
  \div\,\sigmat(\vt,p) = 0,\quad \div\,\vt = 0\text{ in } \FL(t),
\end{equation}
where $\FL(t)$ is the (moving) fluid domain, the  lumen,
$\rho_f\approx \unit[1.06]{g cm^{-3}}$ the density of the blood and
$\sigmat(\vt,p) = \rho_f\nu_f (\nabla\vt+\nabla\vt^T) - pI$
the Cauchy stress tensor, depending on velocity $\vt$ and pressure $p$,
with the kinematic viscosity $\nu_f \approx \unit[0.03]{cm^2
  s^{-1}}$. The vessel walls 
are governed by an elastic material 
\begin{equation}\label{solidproblem:0}
  J\rho_s \partial_t \vt - \div \Big(\Ft\Sigmat \Big) = 0,\quad
  \vt = \partial_t \ut\text{ in }\SO,
\end{equation}
where $\rho_s$ is the fluid's density (in current configuration),
$\vt$ the velocity, $\ut$ the deformation, $\Ft:=I+\nabla\ut$ the
deformation gradient with determinant $J:=\det\,\Ft$. By $\SO$ we
denote the Lagrangian reference configuration. By $\Sigmat$ we denote
the Piola Kirchhoff stresses. The proper modeling of the stresses
within vessel walls is under active research~\cite{Holzapfel2000}. In
particular 
there is still little knowledge on the degree of complexity that is
required for accurately predicting the behavior of the coupled
system.
To incorporate growth of the stenosis in the context of
fluid-structure interactions, the technique of a
multiplicative decomposition of the deformation gradient
\[
\Ft = \Ft_e \Ft_g(c),\quad \Ft_e = \Ft\Ft_g(c)^{-1},\quad
V_0 \xrightarrow{\Ft_g} V_g \xrightarrow{\Ft_e} V
\]
into active deformation  $\Ft_g(c)$ and elastic response 
$\Ft_e$ can be applied, see~\cite{RodriguezHogerMcCulloch1994,FreiRichterWick2016}. The idea is to introduce an
intermediate configuration that includes the growth $\hat S\to\hat
S_g$ and that is mediated by $\Ft_g(c)$ depending directly on the
exterior growth trigger $c$ but that is not physical, i.e. it is
stress free but not necessarily free of strain. The stresses then
depend on the elastic part only, to be precise on $\Ft_e =
\Ft_g(c)^{-1}(I+\nabla\ut)$. Such models are successfully used in
describing the formation of
plaques~\cite{YangJaegerNeussRaduRichter2015,YangRichterJaegerNeussRadu2017}.  

In this work we considerably simplify the model by neglecting all
elastic effects. The Navier-Stokes equations are solved in the domain
$\FL$ that directly depends on a growth variable $c$ by prescribing
normal growth 
\[
\partial\FL\big(c(t)\big) = \{ \xt-c(\xt,t)\cdot \nt_{\hat
  \FL}(\xt)\,:\, \xt\in\partial \hat\FL\}, 
\]
where $\hat\FL$ is the non-grown fluid domain in reference state and
$\nt_{\hat\FL}$ the outward facing unit normal vector. The
description of the coupled problem we will be based on an ALE
formulation, where all quantities are given on the undeformed
reference domain $\hat\FL$, see~\cite[Chapter 5]{Richter2017}.
This reference domain is a
straight pipe of length $\unit[7]{cm}$ and diameter
$\unit[0.2]{cm}$. A typical curvature, irregularities, the effect of
the stent and in particular of the stenosis will be augmented by the ALE
deformation $T(t):\hat\FL\to\FL(t)$.

The growth variable $c$ will live 
on the surface $\partial\hat\FL$. The evolution of $c$ is governed by
a simple surface diffusion equation
\begin{equation}\label{growth:0}
  d_t c  - \lambda_c \Delta_\Gamma c = 
  R(c,\sigmat)
  \text{ on }\partial\hat\FL,
\end{equation}
with the Laplace Beltrami operator $\Delta_\Gamma$ and a small diffusion
constant $\lambda_c\approx \unit[5\cdot 10^{-7}]{m^2/s}$. Due to the
very slow evolution of the plaque,  the motion of
the evolving surface can be neglected in the temporal
derivative. There is no experimental data on the role of diffusion and
the size of $\lambda_c$. We will hence consider $\lambda_c$ as a
procedure for stabilization and choose is small enough to cancel any
effects on the macroscopic evolution of the growth. In lack of
relevant parameters equation~(\ref{growth:0}) can be considered to be
dimensionless. By $R(c,\sigmat)$ we denote the coupling term
triggering growth of the stenosis
\begin{equation}\label{reaction}
  R(c,\sigmat;\xt)= \frac{\alpha}{1+\beta c(\xt)}
  \gamma(\sigma_{WSS}\big(\sigmat(\xt);\xt)\big).
\end{equation}
The parameter $\alpha$ controls the rate of the stenosis growth and it
can be considered as the scale parameter separating the fast scale of
the fluid problem from the slow scale of the growth, by $\beta$ we
control some saturation of the growth. By $\sigma_{WSS}$ we denote the
wall shear stress that is acting close to the tips of the stent at
$s_0$ and $s_1$ (in direction of the main flow direction $\xt_1$, where  
injuring of the vessel wall will trigger stenosis growth
\[
\sigma_{WSS}(\sigmat;\xt) = |\sigmat(\xt)\nt(\xt)\cdot\nt(\xt)|
\Big(\Theta(s_0;\xt_1)+\Theta(s_1;\xt_2)\Big),
\]
with
\[
\Theta(s;x) = \Big(1+\exp\big(2(s_0-1-x)\big)\Big)^{-1}
\Big(1+\exp\big(2(x-s_0-1)\big)\Big)^{-1}.
\]
Only wall shear stresses in a certain range above and below activation
limits are responsible for plaque growth, hence we introduce the
scaling function $\gamma(\cdot)$ as 
\[
\gamma\big(S\big) =
\Big(1+\exp\big(3(\sigma_{min}-S)\big)\Big)^{-1}
\Big(1+\exp\big(3(S-\sigma_{max})\big)\Big)^{-1}.
\]
%


\subsection{Parameters}

All computations are carried out on
the reference domain, a vessel of diameter $0.2\unit{cm}$ and length
$7\unit{cm}$. Deformations, imposed by the stent $T_{stent}$, the general
curvature of the configuration $T_{geometry}$ and the stenosis
$T_{stenosis}$ are realized by mappings
\[
T = T_{geometry}\circ T_{stenosis}\circ T_{stent}. 
\]
All units are given in $\unit{cm},\unit{g},\unit{s}$.

$T_{stent}$ models the impact of the stent, a slight extension of the
vessel at the tips $s_l$ and $s_r$
\begin{equation}\label{stentgeo}
  T_{stent}(x) 
  ={\small\begin{pmatrix}x_1 \\ 0 \\ 0  \end{pmatrix}}
  + \Big(1+\rho_{stent}\mathrm{e}^{-\gamma_{stent}(x_1-s_0)^2}
  +\mathrm{e}^{-\gamma_{stent}(x_1-s_1)^2}\Big)
  {\small \begin{pmatrix}0   \\x_2 \\x_3 \end{pmatrix}}
\end{equation}
with $\rho_{stent}=0.1$ and $\gamma_{stent}=50$. Growth of the
stenosis is assumed to be in normal direction only. We prescribe
$T_{stenosis}$ by the simple relation
\[
T_{stenosis}(c;x) =
{\small\begin{pmatrix}  x_1 \\ 0 \\ 0\end{pmatrix}}
+\big(1-c(x)\big)
{\small  \begin{pmatrix}  0 \\ x_2\\ x_3\end{pmatrix}}.
\]
The overall vessel geometry is curved in the
$x/y$ plane for $x_1<s_m=\unit[3.5]{cm}$ which is the left half of the
vessel and in the $x/z$ plane for $x_1>s_m$
\[
\begin{aligned}
  T_{geo}(x)\Big|_{x_1<s_m}\hspace{-0.3cm}&=
  {\small \begin{pmatrix}
    x_1-\tau(x_1)\big(1+\tau(x_1)^2\big)^{-\frac{1}{2}} x_2 \\
    \tau'(x_1)+\big(1+\tau(x_1)^2\big)^{-\frac{1}{2}} x_2 \\
    x_3
  \end{pmatrix}},\;
  T_{geo}(x)\Big|_{x_1>s_m}\hspace{-0.5cm}&=
  {\small \begin{pmatrix}
    x_1-\tau(x_1)\big(1+\tau(x_1)^2\big)^{-\frac{1}{2}} x_3 \\
    x_2 \\ 
    \tau'(x_1)+\big(1+\tau(x_1)^2\big)^{-\frac{1}{2}} x_3 \\
  \end{pmatrix}}
\end{aligned}
\]
where $\tau(x_1)$ describes the center-line of the deformed vessel,
given by
$
\tau(x_1) = 4\cdot 10^{-3}(x_1-s_m)^4.
$
The mapping is chosen to give a curvature that is realistic in
coronary arteries with a straight middle-section describing the
stented area. As further parameters we consider the fluid density
$\rho_f=\unit[1.06]{g\cdot cm^{-3}}$, the viscosity $\nu =
\unit[0.03]{cm^2\cdot s^{-1}}$. The stent starts at
$s_0=\unit[2]{cm}$, extends over $\unit[3]{cm}$ to
$s_1=\unit[5]{cm}$. The geometric parameters for the impact of the
stent, see~\ref{stentgeo}, are $\gamma_{stent}=50$ and finally, the
reaction term uses the limits $\sigma_{min}=5$ and $\sigma_{max}=8$. 

The flow problem is driven by enforcing a periodic relative pressure
profile (inflow to outflow) condition that is inspired from the usual
pressure drops in stented coronary arteries suffering from a
stenosis. On the inflow boundary $\Gamma_{in}$ we prescribe the
time-periodic average pressure
\[
P_{in}(t) = \begin{cases}
  10+25 t& 0\le t<\unit[0.4]{s}\\
  140/3 - 200t/3 & \unit[0.4]{s}\le t <\unit[0.7]{s}\\
  100/3t -70t/3& \unit[0.7]{s}\le t <\unit[1]{s}
\end{cases},\quad \text{periodically extended over }[0,1]
\]

\subsection{ALE Formulation and Discretization}

Based on the mapping $T(x)=T_{geometry}(x)\circ T_{stenosis}(x)\circ
T_{stent}(x)$ the Navier-Stokes equations and the surface growth
equation are transformed to ALE coordinate, e.g. by introducing
reference values $\hat v(\hat x,t) = v(x,t)$, $\hat p(\hat
x,t)=p(x,t)$ and $\hat c(\hat x,t) = c(x,t)$. The resulting set of
equations is given on the reference domain $\hat\FL$ and in
variational formulation it takes the form
\begin{equation}\label{ALE}
  \begin{aligned}
    \Big(J\rho_f \big(\partial_t \hat\vt +
    (\hat
    \Ft^{-1}\hat\vt\cdot\hat\nabla)\hat\vt\big),\phi\Big)_{\hat\FL}  
    +\Big( J\hat\sigmat\Ft^{-T},\hat\nabla\hat\phi\Big)_{\hat\FL}& = 0
    \\
    \Big(J \Ft^{-1}:\hat\nabla\hat\vt,\xi\Big)_{\hat\FL}=0,\qquad
    \Big(c',\psi\Big)_{\hat\partial\FL} + \Big(\lambda_c \nabla_\Gamma
    c,\nabla_\Gamma \psi\Big)_{\hat\partial\FL} & = R(\hat
    c,\hat\sigmat).  
  \end{aligned}
\end{equation}
Several simplifications in comparison to an exact ALE formulation have
been applied: due to the very slow evolution of the surface we neglect
inertia terms by its motion. Further, since surface diffusion will
only serve as numerical stabilization we refrain from an exact
transformation of the surface Laplace. 

The discretization of system~(\ref{ALE}) is by standard techniques. In
time, we use the $\theta$-time stepping method
\[
u'=f(t,u)\quad\rightarrow\quad
u_n-u_{n-1} = \Delta t \theta f(t_n,u_n)+\Delta
t(1-\theta)f(t_{n-1},u_{n-1}), 
\]
with constant step sizes $\Delta$ and the choice
$\theta=\frac{1}{2}+{\cal O}(k)$ to achieve second order accuracy with
good stability properties,
see~\cite{LuskinRannacher1982,RichterWick2015_time}. Spatial
discretization is by means of stabilized equal order tri-quadratic
finite elements on a hexahedral mesh. For stabilization of the inf-sup
condition and of convective regimes the local projection stabilization
is used~\cite{BeckerBraack2001,BeckerBraack2004}. The surface pde is
continued into the fluid domain and can be considered as a weakly
imposed boundary condition. We refer to~\cite{Richter2017} for details
on the discretization and implementation in Gascoigne
3D~\cite{Gascoigne3D}. 

\section{Temporal multiscales}\label{sec:tempmulti}

The big challenge of the candy wrapper problem is in the range of
temporal scales that must be bridged. While the flow problem is driven
by a periodic flow pattern with period $\unit[1]{s}$ the growth of the
stenosis takes months. The growth model comprises the parameter
$\alpha$, see~(\ref{reaction}) that indicates exactly this scale separation,
since $|R(c,\sigmat)| = {\cal O}(\alpha)$. In~\cite{FreiRichter2019}
we have recently 
introduced and analysed a temporal multiscale scheme for exactly such
long-scale / short-scale problems governed by a pde/ode system and
driven by a periodic-in-time micro process. Here we extend this
technique for handling 3d pde/pde couplings. 

We briefly sketch the layout of the multiscale approximation. To begin
with, we identify the growth parameter $c(\xt,t)$ as the main variable of
interest. Furthermore, as we are interested in the long term behavior
of the growth only, we introduce the (locally) averaged growth
variable 
\begin{equation}\label{averaged:c}
  \bar c(\xt,t) = \int_{t}^{t+\unit[1]{s}} c(\xt,s)\,\text{d}s,
\end{equation}
where the averaging extends over one period only.

Next, to decouple slow and fast scales we make the essential
assumption that the flow problem on a fixed domain $\FL(\bar c_{f})$,
where $\bar c_f:=\bar c(t_f)$ for one point in time $t_f$ 
admits a periodic in time solution
\begin{multline}\label{per}
    \Big(J(\bar c_f)\rho_f \big(\partial_t^{\bar c_f} \vt^{\bar c_f} +
    (\Ft(\bar c_f)^{-1}\vt^{\bar c_f}\cdot\nabla)\vt^{\bar
      c_f}\big),\phi\Big)_{\FL}   \\
    +\Big( J(\bar c_f)\sigmat^{\bar c_f}\Ft(\bar c_f)^{-T},\nabla\phi\Big)_{\FL}
    +    \Big(J(\bar c_f) \Ft(\bar c_f)^{-1}:\nabla\vt^{\bar
      c_f},\xi\Big)_{\FL}  = 0\\
    \vt^{\bar c_f}(\cdot,0) =
    \vt^{\bar c_f}(\cdot,1) \qquad\qquad\qquad\qquad
\end{multline}
Only very few theoretical results exist on periodic solutions to the
Navier-Stokes equations, see~\cite{GaldiKyed2016}. They only hold in
the case of 
small data which is not given in the typical candy wrapper
configurations with Reynolds numbers going up to about
$Re=1\,000$. Computational experiments however do suggest the 
existence of stable periodic solutions in the regime of interest.

\textbf{Multiscale Algorithm}
Given such periodic solutions, the computational multiscale method is
based on a subdivision of $I=[0,T]$ (where $T\approx \unit{months}$ is
 large) into macro time-steps $t_n$ for $n=0,\dots,N$ with $t_0=0$ and
 $T_N=T$ and the step size $K=t_n-t_{n-1}$. 
The small interval of periodicity $I_P=[0,1]$
is partitioned into micro time-steps $\tau_n$ for $n=0,\dots,M$ with
$\tau_0=0$, $\tau_M=1$ and the step size $k=\tau_m-\tau_{m-1}\ll K$.
A simple explicit/implicit multiscale iteration is then as follows:
\begin{algo}[First order explicit/implicit multiscale iteration]
  Let  $\bar c_0$ be the initial value for the slow component. For
  $n=1,2,\dots$ iterate
  \begin{enumerate}
  \item Solve the periodic flow problem $(\vt^{\bar
    c_{n-1}},p^{\bar c_n})$ on the domain $\FL(\bar c_{n-1})$
  \item Compute the average of the reaction term
    \begin{multline*}
      \bar R(\bar c_{n-1}):=
      \int_0^1 R(\bar c_{n-1},\sigmat^{\bar c_{n-1}}(s);\xt)\,\text{d}s\\
      = \frac{\alpha}{1+\beta \bar c_{n-1}(\xt)}
      \int_0^1 \gamma\Big(\sigma_{WSS}\big(\sigmat^{\bar
        c_{n-1}}(\xt,s);\xt\big)\Big)\,\text{d}s 
    \end{multline*}
  \item Make an semi-explicit step of the stenosis growth problem
    \[
    K^{-1}\big(\bar c_n-c_{n-1},\psi\big)_{\partial\hat\FL}
    +\big(\lambda_c\nabla_\Gamma \bar
    c_{n},\nabla_{\Gamma}\psi\big)_{\partial\hat\FL} 
    =  \big(\bar R(\bar c_{n-1}),\psi\big)_{\partial\hat\FL}
    \]
  \end{enumerate}
\end{algo}
The discretization of the growth problem in Step 3. can easily be
replaced by a second order explicit scheme like the Adams-Bashforth
formula, see~\cite{FreiRichter2019}. A fully implicit time-integration
can be realized by adding a sub-iteration for steps 2-4. However,
since the diffusion parameter is very small, explicit schemes are
appropriate in this setting. 

Within every step of the iteration it is necessary to solve the
periodic-in-time flow problem (even multiple solutions are required in
a fully implicit setting). This is the main effort of the resulting
scheme, since the sub interval $[0,1]$ must be integrated several
times to obtain a suitable periodic solution. In principle it is
possible to just compute several cycles of the periodic problem until
the periodicity error
\[
\|\vt^{\bar c_n}(T+\unit[1]{s})-\vt^{\bar c_n}(T)\|< \epsilon_P
\]
falls below a given threshold $\epsilon_P>0$. Usually however this
error is decreasing with an exponential rate only that depends on
parameters like the viscosity and the domain size. For
acceleration several methods are discussed in literature, based on
optimization problem~\cite{RichterWollner2018}, on the idea of the
shooting method~\cite{JianBieglerFox2003}, on
Newton~\cite{HanteMommerPotschka2015} or on space time
techniques~\cite{PlatteKuzminFredebeulTurek2005}.  
Here we quickly present a very efficient novel scheme that converges
with a fixed rate that does not depend on any further parameters. We
note however that although the computational efficiency is striking,
the theoretical validation extends to the linear Stokes equation only,
see~\cite{Richter2019}.

\textbf{Solution of the periodic flow problem}
The idea of the averaging scheme for the rapid identification of
periodic flow problems is to split the periodic solution into
average and oscillation, see also~\cite{Richter2019}
\[
\vt^\pi(t) = \bar \vt^\pi + \tilde \vt^\pi(t),\quad \int_0^1
\tilde \vt^\pi(s)\,\text{d}s=0. 
\]
In a nonlinear problem like the Navier-Stokes equations it is not
possible to separate the average from the oscillations. But, by
averaging the Navier-Stokes equation, we derive 
\[
-\div\,\bar\sigmat_f^\pi + (\bar\vt^\pi\cdot\nabla)\bar\vt^\pi
= \underbrace{-\int_0^1 \Big\{ (\tilde\vt^\pi(s)\cdot\nabla) \bar\vt^\pi
+
(\bar\vt^\pi\cdot\nabla)\tilde\vt^\pi(s)\Big\}\text{d}s}_{=:N(\bar\vt^\pi,\tilde\vt^\pi)},\quad 
\div\,\bar\vt^\pi=0. 
\]
If we average a solution $(\vt(t),p(t))$ to the Navier-Stokes problem
for arbitrary initial $\vt_0$ (that does not yield the periodic
solution) we get
\[
-\div\,\bar\sigmat_f + (\bar\vt\cdot\nabla)\bar\vt
= \vt(0)-\vt(1) + 
N(\bar\vt,\tilde\vt),\quad 
\div\,\bar\vt=0. 
\]
The  difference $\wt:=\vt^\pi-\vt$, $q:=p^\pi-p$ between dynamic 
solution and periodic solution satisfies the averaged equation
\begin{multline*}
  -\div\,\bar\sigmat_f(\bar\wt,\bar q) + (\bar\wt\cdot\nabla)\bar\wt
  +(\wt\cdot\nabla)\bar\vt + (\bar\vt\cdot\nabla) \wt\\
  = \vt(1)-\vt(0) +
  N(\bar\vt^\pi,\tilde\vt^\pi)
  -
  N(\bar\vt^\pi,\tilde\vt^\pi),\quad 
  \div\,\bar\wt=0. 
\end{multline*}
We assume that we start with a good guess $\vt$ that is already close
to the periodic solution $\vt^\pi$, i.e. $\|\wt\|$ is small. If no
initial is available, e.g. in the very first step of the multiscale
scheme, we still can perform a couple for forward simulations. Given
that $\|\wt\|$ is small, we will neglect both the nonlinearity
$(\bar\wt\cdot\nabla)\bar\wt$ and all fluctuation terms
$N(\cdot,\cdot)$ involving the oscillatory parts. We approximate
the difference between average of the dynamic solution and average of
the desired periodic solution by the linear equation
\begin{equation}\label{avg}
  (\wt\cdot\nabla)\bar\vt + (\bar\vt\cdot\nabla) \wt
  -\div\,\bar\sigmat_f(\bar\wt,\bar q)
  = \vt(1)-\vt(0)
\end{equation}
The averaging scheme for finding the periodic solution is then given
by the following iteration.

\begin{algo}[Averaging scheme for periodic-in-time problems]
  Let $\vt^0_0$ be a guess for the initial value. If no approximation
  is available, $\vt^0_0$ can be obtained by computing several cycles
  of the dynamic flow problem. For $l=1,2,\dots$ iterate
  \begin{enumerate}
  \item Based on the initial $\vt^l(0)=\vt^{l-1}_0$ solve once cycle of
    the dynamic flow problem on $I_P=[0,1]$.
  \item Solve the averaging equation for $\bar\wt^l$ and  $\bar q^l$,
    equation~(\ref{avg})
  \item Update the initial value by correcting the average
    \[
    \vt^l_0:=\vt^l(1) + \bar\wt^l.
    \]
  \end{enumerate}
\end{algo}
The analysis of this averaging scheme is open for the Navier-Stokes
equations but simple for linear problems with symmetric positive
definite operator like the Stokes equations. Here the convergence
estimate
\[
\|\vt^l_0- \vt^\pi_0\|\le \rho_{avg} \cdot \|\vt^{l-1}_0-\vt^\pi_0\|
\]
holds, with $\rho_{avg}<0.3$ in the continuous and $\rho_{avg}<0.42$
in the discrete setting, for further results we refer
to~\cite{Richter2019}. 

\section{Numerical results}

We present a numerical study on the multiscale scheme and give a first
discussion on its accuracy and efficiency. In~\cite{FreiRichter2019}
simpler two dimensional problems have been studied that also allow for
resolved simulations such that a direct comparison of computational
times for forward simulations and multiscale simulations can be
performed. These demonstrated speedups reaching from $1:200$ to
$1:10\,000$. Here it was shown that the multiscale scheme benefits
from larger scale separation. To be precise: to reach the same
relative accuracy in a multiscale computation as compared to a direct
forward computation, the speedup behaves like $1:\alpha^{-1}$.

Before presenting results for the multiscale method we briefly discuss
the averaging scheme for finding periodic solutions 

\subsection{Convergence of the averaging scheme for periodic
  flow problems}

\begin{table}[t]
  \begin{center}
    \begin{tabular}{lcccc|cccc}
      \toprule
      $\nu$ & 0.1 & 0.05& 0.025 &&& 0.1 & 0.05& 0.025 \\
      &\multicolumn{3}{c}{forward}&&
      &\multicolumn{3}{c}{averaging}\\
      \midrule
      cycles & 40 & 74 & 140 &&  & 15 & 15 & 18\\
      \bottomrule
    \end{tabular}
  \end{center}
  \caption{Number of cycles required to reduce the periodicity error
    to $\|\vt(t+\unit[1]{s})-\vt(t)\|<10^{-8}$ for the direct forward
    simulation and the averaging scheme. Variation in the viscosity
    $\nu$.}
  \label{tab:periodic}
\end{table}

We consider a 3d problem that is inspired by the driven cavity
problem. On the cube $\Omega=(-2,2)^3$ we drive the Navier-Stokes
equation by a $1$-periodic forcing
\[
\ft(\xt,t) = \frac{\sin\big(2\pi t\big)}{6}
\begin{pmatrix}
  3\tanh(\xt_2)\\
  2\tanh(\xt_3)\\
  \tanh(\xt_1)\\
\end{pmatrix}
\]
Since the data is periodic in time we can expect to obtain a
time-periodic solution (if the Reynolds number is sufficiently
small). In Tab.~\ref{tab:periodic} we show the performance of the
averaging scheme in comparison to a simple forward iteration. We give
the number of cycles required to reach the periodicity error 
 $\|\vt(T_n+\unit[1]{s}) - \vt(T_n)\|<10^{-8}$. The results show a
strong superiority of the averaging scheme, both in terms of
robustness (with respect to $\nu$) and in terms of the overall
computational complexity.  For the forward iteration, the number of
cycles approximately doubles with each reduction of $\nu$. The
performance of the averaging scheme slightly deteriorates for 
$\nu=0.025$ due to the higher Reynolds number regime. For $\nu=0.01$
we cannot identify a stable periodic solution. The computational
overhead of the averaging scheme is very low, one additional
stationary problem must be solved in each cycle. A detailed study of
the 
averaging scheme with an analysis of the sensitivity to various
further parameters is given in~\cite{Richter2019}. 

\subsection{Simulation of the candy wrapper problem}

\begin{figure}[t]
  \begin{center}
    \includegraphics[width=\textwidth]{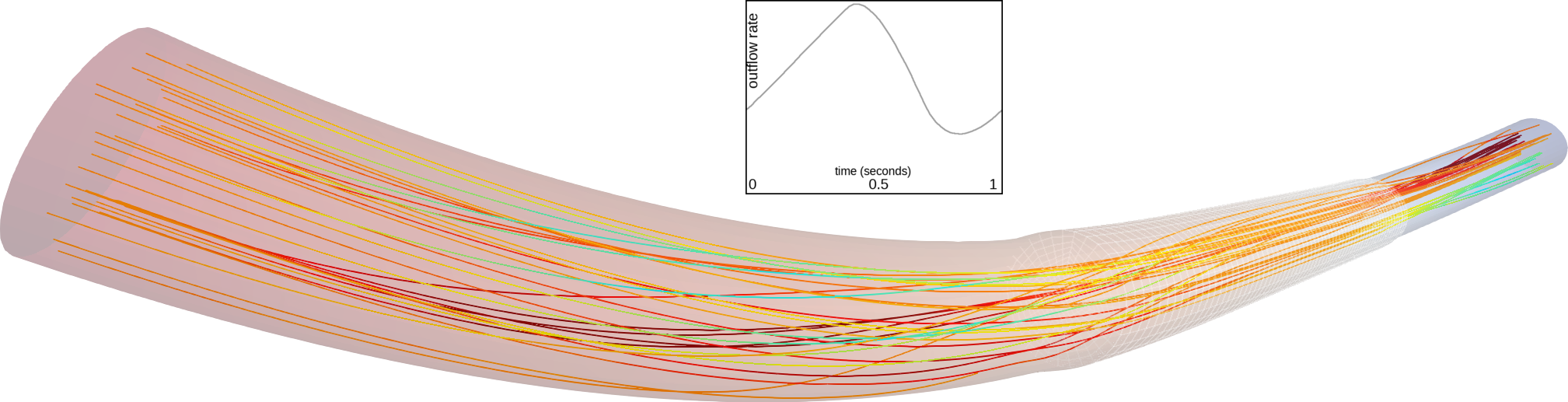}
    \includegraphics[width=\textwidth]{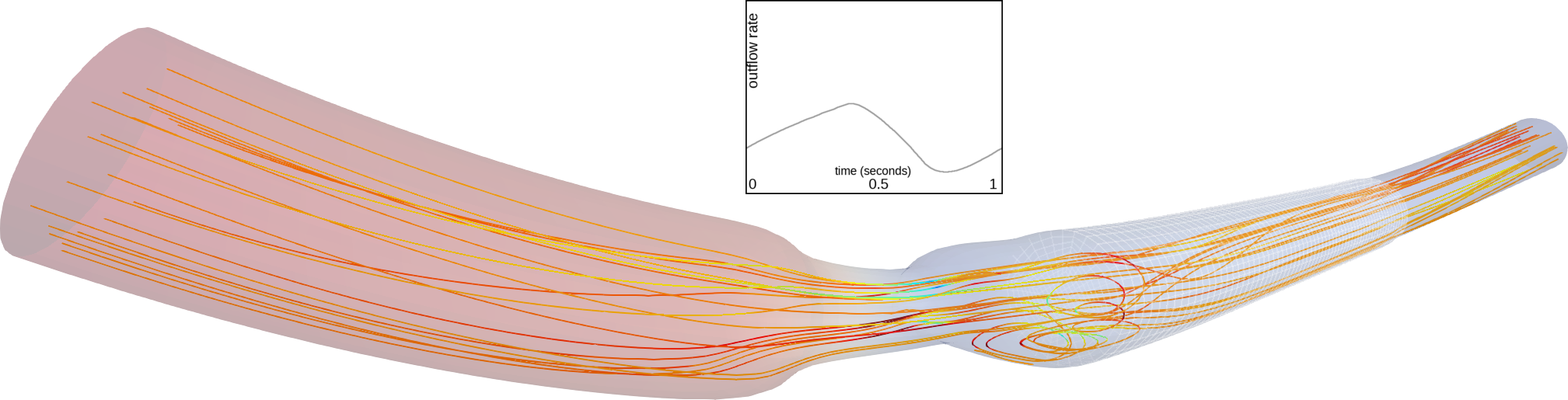}
    \includegraphics[width=\textwidth]{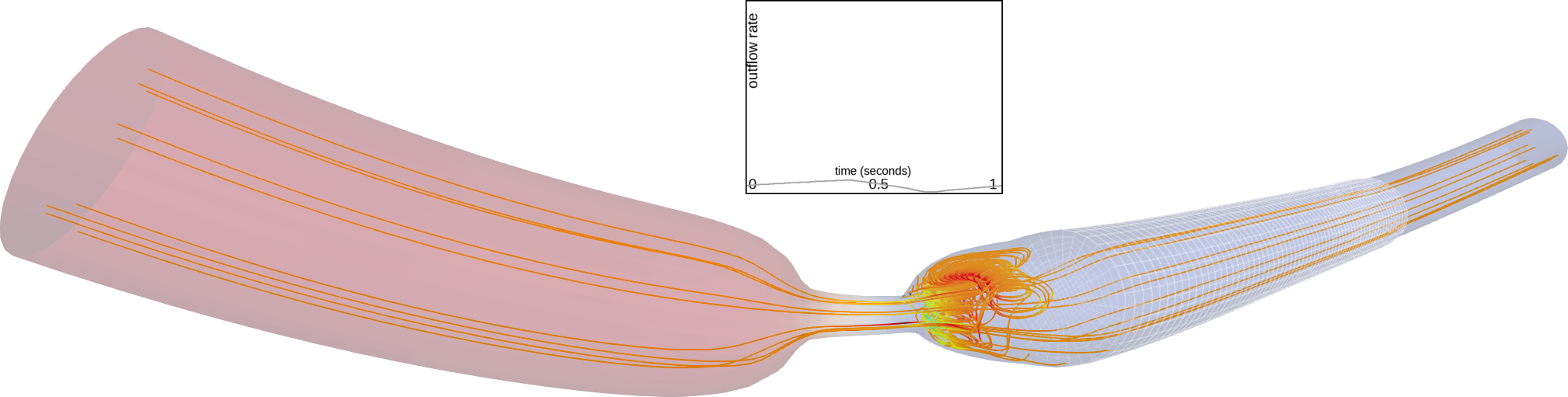}
  \end{center}
  \caption{Development of the stenosis at initial time, 
    at $T=\unit[33]{days}$ and $T=\unit[67]{days}$. The average
    and the oscillation of the flow rate get smaller while the
    stenosis develops. }
  \label{fig:candy}
\end{figure}

Fig.~\ref{fig:candy} shows the evolution of the stenosis at three
different points in time. In addition we show the outflow rate as
function over time (one period). Several effects known from the
medical practice can be identified: The growth of the stenosis is
non-symmetric and mostly centered on the inflow-tip of the stent. This
shows the necessity of considering full three dimensional
models. Further, the simulations show an extension and growth of the
stenosis to both sides which is also typical. Since the flow is
pressure driven, the outflow rate decreases with the development of
the stent. 

\begin{table}[t]
\begin{center}
  \begin{tabular}{rrcrc|crrcr}
    \toprule
    $K$ & $k$ & $J_{out}$&time&&&$K$ & $k$ & $J_{out}$&time \\
    \midrule
    144\,000&0.02&0.9359&$\unit[9]{min}$&&& 72\,000&0.04 &0.9132&$\unit[15]{min}$\\
    72\,000&0.02 &0.9138&$\unit[18]{min}$&&& 72\,000&0.02 &0.9138&$\unit[18]{min}$\\
    36\,000&0.02 &0.9043&$\unit[40]{min}$&&&72\,000&0.01
    &0.9140&$\unit[41]{min}$\\ 
    \midrule
    \multicolumn{2}{l}{extra $K\to 0$}&\phantom{X} 0.8971 (1.22) &
    \phantom{X}$\unit[55]{years}^{(*)}$&&&
    \multicolumn{2}{l}{extra $k\to 0$}&\phantom{X} 0.9131 (1.58)\\
    \bottomrule
  \end{tabular}
\end{center}
\caption{Outflow $J_{out}$ at time $T\approx \unit[18]{days}$ and
  extrapolation including numerical convergence order for $K\to 0$
  ($k$ fixed) and $k\to 0$ ($K$ fixed). 55 years $(*)$ computational
  time result from a projection of the 
  computational time for a resolved simulation without the multiscale
  scheme. } 
\label{tab:res}
\end{table}

In Table~\ref{tab:res} we compare the results of the multiscale scheme
for different values of $k$ and $K$. We observe convergence in both
parameters. Numerical extrapolation yields ${\cal
  O}(k^{1.58}+K^{1.58})$, slightly off the expected rates ${\cal
  O}(k^2+K)$. We also indicate the computational times required for
running the multiscale scheme till $T\approx \unit[18]{days}$. A
corresponding resolved simulation would require about
$\unit[55]{years}$ computational time. This value is predicted based
on the average time for computing a complete cycle of the periodic
problem and based on an average three iterations required for
approximating the periodic flow problem. Assuming that the
extrapolated value for $K\to 0$ is accurate, the simulation based on
$K\approx \unit[36\,000]{s}$ carries a multiscale error of about
$1\%$. This approximation is achieved in $\unit[40]{min}$ instead of
$\unit[55]{years}$. The results in Table~\ref{tab:res} indicate that
it is worthwhile to consider a second order time stepping scheme for
the plaque growth problem, since the error in $K$ dominates. We refer
to~\cite{FreiRichter2019} for a realization in the context of a
pde / ode long-scale / short scale problem. 

\section{Outlook and discussion}

We have demonstrated a numerical framework for simulating complex
multiphysics / multiscale problems in hemodynamics. For the first time
we could demonstrate an efficient numerical scheme for a  long-scale /
short-scale problem coupling different partial differential
equations. We are able to include both temporal and spatial effects in
bio-medical growth applications. The combination of a temporal
multiscale method with fast solvers and efficient discretizations for
the (periodic) micro problems gives substantial speedups such that
three dimensional problems can be treated. Two main challenges remain
for future work:

\textbf{Fluid-structure interactions}
The main challenge in including elastic vessel walls lies in the
increased complexity of the resulting  system due to nonlinearities
coming from the domain motion and the coupling to the hyperbolic solid
equation that, by introducing the deformation as additional variable,
blows up the problem size. In hemodynamical applications the coupling
is governed by the added mass
instability that usually calls for strongly coupled solution
approaches,
see~\cite{CausinGerbeauNobile2005,HeilHazelBoyle2008}. Although some
progress has been made in recent
years~\cite{Richter2015,AulisaBnaBornia2018,JodlbauerLanger,FailerRichter2019},
the design of efficient solvers for the resulting algebraic problems is still not
satisfactory.

Considering monolithic solution approaches in combination with
Newton-Krylov solvers make the use of large time steps possible. In
all of the just mentioned approaches for designing linear solvers it
has shown to be essential to partition the linear system when it
actually comes to inversion of matrices, either within a
preconditioner or within a multigrid smoother. This is mainly due to
the very large condition numbers of the coupled system matrix that by
far exceeds those of the subproblems,
see~\cite{Richter2015,AulisaBnaBornia2018}.

A second difficulty coming with fluid-structure interactions lies in
the derivation of the effective growth equation described in
Section~\ref{sec:tempmulti}. If elastic fluid-structure interactions
are taken into account, the domain undergoes oscillations in the scale
of the fast problem, i.e. during each pulsation of the blood
flow. However, we can nevertheless introduce the averaged growth
variable $\bar c(\xt,t)$ as in~(\ref{averaged:c}) and simply average
the growth equation, the this equation of Eq.~(\ref{ALE}) as this is
stated on the fixed reference domain. We note however that we have
chosen a very simple growth model given as surface
equation. Considering the detailed system introduced
in~\cite{YangJaegerNeussRaduRichter2015}, growth takes place within
the solid, which is a three dimensional domain
$\SO(t)\subset\mathds{R}^3$ undergoing deformation from the coupled
fsi problem. A corresponding equation mapped to the fixed reference
domain $\hat\SO$ (taken from
\cite{YangJaegerNeussRaduRichter2015}) reads 
\[
\Big(
\frac{\partial}{\partial t}(J\hat c),\psi\Big)_{\hat\SO}
+\Big(\lambda_c J\Ft^{-1}\Ft^{-T}\nabla \hat c,\nabla\psi
\Big)_{\hat\SO} = R(\hat c,\hat\sigmat). 
\]
Since $J$ and $\Ft$ oscillate with the frequency of the fast scale
problem, derivation of an effective equation is still subject  to
future work. 

\textbf{Patient specific simulation}
The second open problem is to incorporate patient specific data into
the simulations for generating specific predictions. Flow and geometry
data can easily be measured during the stenting process. This process
however is strongly invasive and causes subsequent adaptions of the
vessel and the surrounding tissue interacting with the stent. Further
data on the resulting configurations are not easily available without
additional interventions. With a diameter of only a few millimeters, coronary
arteries are small, such that measurements at good accuracy cannot be
obtained.

\textbf{Medical application}
The edge stenosis accompanying the implementation of DES (Drug Eluting Stents) is great starting point for development of the further numerical experiments in the field of the plaque formation and biochemical processes ongoing in the vessel walls exposed to other types of interventions.  
Explosive growth of the intravascular interventions in recent decade is, inevitably, going to demand more advanced studies on the nature of vascular wall response to the implantable devices~\cite{Mack2019}. Novel numerical methods may also shade new light on well-established surgical procedures and augment the awareness of the potential benefits or hazards that are not yet fully understood or identified~\cite{Thourani}. On the other hand the population of the patients is changing dramatically and that process is soon to accelerate. According to the recent report published by European Commission, diseases of the circulatory system are the most common cause of death in elderly population aged over 75 years~\cite{Eurostat}. In addition to that gruesome information the ageing of the European population in the years to come is growing concern of the governments. Poland belongs to the group of the countries that may become affected by the population ageing the most~\cite{Giannakouris}. Due to that we face the necessity of development the most efficient treatment strategies for the elderly population. One of those treatment procedures is TAVR (Transcatheter Aortic Valve Replacement). The procedure addresses aortic valve stenosis that is quite often ailment in the aforementioned group of patients. By application of the fluid structure interaction methods it might be possible to tailor the design of the medical devices to the stiffer tissues usually present in the elderly patients in the way that may augment long time outcome of the procedure. Just such a small improvement may diminish the risk of repeated procedures undertaken in frail patients. 

The methodology presented in our work should also find it’s application in optimization of the classic surgery for the coronary artery disease. The position of the vascular anastomosis in relation to the existing vascular wall lesions may find new rationale when understood through the knowledge of the mechanotranduscion phenomena. Also the strategic planning of the target vessels and ``landing sites'' for the aorto-coronary by-pass grafts may find its' new understanding. Those perspective studies could be undertaken only by the means of model based planning.

\section*{Acknowledgement}

We acknowledge support by the Federal Ministry of Education and
Research of Germany (project number 05M16NMA). 


\end{document}